
\input amstex
\magnification=1200

\documentstyle{amsppt}
\pagewidth{5.5 in} \pageheight{7.5 in}

\NoRunningHeads
\nologo

\def\R{\Bbb{R}}

\def\C{\Bbb{C}}

\def\End{\operatorname{End}}
\def\supp{\operatorname{Supp}}
\def\spec{\operatorname{Spec}}
\def\Min{\operatorname{Min}}

\def\m{{\frak m}}
\def\n{{\frak n}}
\def\fa{{\frak a}}

\def\iff{{\Longleftrightarrow}}

\def\V{\operatorname{V}}

\def\Spec{\operatorname{Spec}}

\def\Ker{\operatorname{Ker}}
\def\ann{\operatorname{Ann}}
\def\im{\operatorname{Im}}

\def\Coker{\operatorname{Coker}}

\def\dim{\operatorname{dim}}

\def\Hom{\operatorname{Hom}}
\def\Tor{\operatorname{Tor}}
\def\Ext{\operatorname{Ext}}

\def\rank{\operatorname{rank}}
\def\depth{\operatorname{depth}}

\def\HH{\operatorname{H}}

\def\vf{\varphi}
\def\fm{\frak m}
\def\fn{\frak n}
\def\ol{\overline}

\def\Supp{\operatorname{Supp}}

\def\pd{\operatorname{pd}}

\def\V{\operatorname{V}}

\def\Rh{R^{\text{h}}}

\topmatter

\title Ascent of module structures, vanishing of Ext, and extended modules\endtitle
\rightheadtext{}

\author Anders J.~Frankild\footnote{This work was
completed after the untimely death of Anders J. Frankild on 10 June 2007.}
and Sean Sather-Wagstaff
and  Roger Wiegand\footnote{Wiegand's research was partially supported by Grant
04G-080 from the National Security Agency.}
\endauthor

\address Anders J.\ Frankild, University of Copenhagen, Institute for Mathematical
Sciences, Department of Mathematics,
Universitetsparken 5, 2100 K\o benhavn, Denmark
\endaddress
\urladdr http://www.math.ku.dk/\~{}frankild/ \endurladdr

\address Sean Sather-Wagstaff,
Department of Mathematical Sciences, Kent State University,
Mathematics and Computer Science Building, Summit Street, Kent OH
44242, USA \endaddress
\curraddr 300 Minard Hall North Dakota State University Fargo ND 58105-5075, USA
\endcurraddr
\email Sean.Sather-Wagstaff\@ndsu.edu \endemail
\urladdr http://math.ndsu.nodak.edu/faculty/ssatherw/ \endurladdr

\address Roger Wiegand, Department of Mathematics, University of Nebraska,
203 Avery Hall, Lincoln, NE 68588-0130, USA \endaddress
\email rwiegand\@math.unl.edu \endemail
\urladdr http://www.math.unl.edu/\~{}rwiegand/ \endurladdr

\date November 20, 2007 \enddate

\abstract
Let $(R,\m)$ and $(S,\n)$ be commutative Noetherian local rings, and let
$\varphi:R\to S$ be a flat local homomorphism such that $\m S = \n$ and
the induced map on residue fields $R/\m \to S/\n$ is an isomorphism.
Given a finitely generated $R$-module $M$, we show that $M$ has an
$S$-module structure compatible with the given $R$-module structure if and
only if $\Ext^i_R(S,M)=0$  for each $i\ge 1$.

We say that an $S$-module $N$ is {\it extended} if there is a finitely
generated $R$-module $M$ such that $N\cong S\otimes_RM$.
Given a short exact sequence $0 \to N_1\to N \to N_2\to 0$ of finitely
generated
$S$-modules, with two of the three modules $N_1,N,N_2$ extended, we obtain
conditions forcing the third module to be extended.  We show that every
finitely generated module over the Henselization of $R$ is a direct summand
of an extended module, but that the analogous result fails for the
$\m$-adic completion.
\endabstract

\dedicatory This paper is dedicated to Melvin Hochster on the
occasion of his sixty-fifth birthday. \enddedicatory

\endtopmatter

\document

\head Introduction \endhead

Suppose $(R,\m)$ and $(S,\n)$ are commutative
Noetherian local rings and $\varphi\colon  R \to S$ is a flat local
homomorphism with the property that the induced homomorphism
$R/\m \to  S/\m S$ is bijective. We consider natural questions of
ascent and descent of modules between $R$ and $S$:
(1) Given a finitely generated $R$-module $M$, when does
$M$ have an $S$-module structure that is compatible with the
$R$-module structure via $\varphi$? (2) Given a finitely generated
$S$-module $N$,  is there a finitely generated $R$-module $M$
such that $N$ is $S$-isomorphic to $S\otimes_RM$, or (3)
$S$-isomorphic to a direct summand of $S\otimes_RM$?

In Section 1 we make some general observations about homomorphisms
$R\to S$ satisfying the condition $R/\m=S/\m S$.  We show that if a
compatible $S$-module structure exists, then it arises in an obvious
way: The natural map $M \to S\otimes_RM$ is an isomorphism. (One
example to keep in mind is that of a finite-length module $M$ when
$S = \widehat R$, the $\m$-adic completion.) Moreover, if $R\to S$
is flat, then $M$ has a compatible $S$-module structure if and only
if $S\otimes_RM$ is finitely generated as an $R$-module.

In Section 2 we prove, assuming that $R\to S$ is flat, that $M$ has
a compatible $S$-module structure if and only if  $\Ext^i_R(S,M)$ is
finitely generated as an $R$-module for $i=1,\dots,\dim_R(M)$. We
were motivated to investigate this implication because of the
following result of Buchweitz and Flenner~\cite{BF} and  Frankild
and Sather-Wagstaff~\cite{FSWa}: A finitely generated $R$-module $M$
is $\m$-adically complete if and only if $\Ext^i_R(\widehat{R},M)=0$
for all $i\ge 1$. Theorem 2.5 summarizes the main results of the
first two sections. Note that it subsumes the result of Buchweitz,
Flenner,  Frankild and Sather-Wagstaff, but our proof is quite
different.

In Section 3 we address questions (2) and (3) and show that (3)
always has an affirmative answer when $S$ is the Henselization, but
not necessarily when $S$ is the $\m$-adic completion.

\head 1. Ascent of module structures \endhead

 Throughout this section $(R,\m)$ and $(S,\n)$ are Noetherian local rings
 and $\varphi\colon R\to S$ is a local ring homomorphism.
 We consider the following condition on $\varphi$:
\smallskip
$(\dagger)$ (i) $\m S = \n$; and (ii) $\varphi(R) + \n = S$, that
is, $\varphi$ induces an isomorphism on residue fields.
\smallskip
\noindent This is equivalent to the following:
The induced homomorphism $R/\m \to S/\m S$ is bijective.
induces an isomorphism on residue fields.)


Familiar examples include the $\m$-adic completion $R \to \widehat R$, the
Henselization $R \to \Rh$, and the natural map
$R \twoheadrightarrow S = R/I$, when $I$ is a proper ideal of $R$.

From (i), it follows immediately that $\m^tS = \n^t$ for all $t$.  Similarly,
the next result shows that
(ii) carries over to powers (though here we need {\it both} (i) and (ii), as
is shown by the example $\C[[T^2,T^3]]\subseteq \C[[T]]$).

\proclaim{Lemma 1.1} If  $\varphi\colon  R \to S$ satisfies
$(\dagger)$, then $\varphi(R) + \n^t = S$ for each $t\ge 1$.
\endproclaim
\demo{Proof} By choosing a composition series, we see that
every  $S$-module of finite length has (the same) finite
length as an $R$-module.  In particular, $S/\n^{t+1}$ has
finite length and therefore is finitely generated as an
$R$-module. We have
$$
\frac{\varphi(R)+\n^t}{\n^{t+1}}+\m\frac{S}{\n^{t+1}} =
\frac{\varphi(R)+\n^t+\m S}{\n^{t+1}} = \frac{\varphi(R) +
\n}{\n^{t+1}} = \frac{S}{\n^{t+1}}.
$$
Nakayama's lemma implies that $(\varphi(R)+\n^t)/\n^{t+1} =
S/\n^{t+1}$. \qed\enddemo

The next result is an indispensable tool for several of our proofs.

\proclaim{Proposition 1.2} Assume $\varphi\colon R \to S$
satisfies $(\dagger)$.
Let $M$ and $N$ be $S$-modules, with $_SN$ finitely generated.  Then
$\Hom_R(M,N) = \Hom_S(M,N)$.
\endproclaim
\demo{Proof} We'll show that $\Hom_R(M,N) \subseteq \Hom_S(M,N)$,
since the reverse inclusion is obvious.  Let $f\in \Hom_R(M,N)$.
Given $x\in M$ and $s\in S$, we want to show that $f(sx) = sf(x)$.
Since $_SN$ is finitely generated, it will suffice to show that
$f(sx)-sf(x) \in \n^tN$ for each $t\ge 1$.

Fix an integer $t\ge 1$, and note the following relations:
$$f(\n^t M)  =
f(\m^tM)  \subseteq \m^tN \subseteq \n^tN$$
    Use (1.1) to choose an
element $r\in R$ such that $\varphi(r)-s\in \n^t$.  Then we have
$$
f(sx) - sf(x) = f(sx) - f(rx) + rf(x) - sf(x) = f((s-\varphi(r))x) +
(\varphi(r)-s)f(x).
$$
It follows that $f(sx) - sf(x)$ is in
$$f((s-\varphi(r))M) + (\varphi(r)-s)N \subseteq f(\n^t M) + \n^tN
= \n^tN. \qed$$
\enddemo

\proclaim{Corollary 1.3}
Let $\vf\colon R\to S$ be a local homomorphism
satisfying $(\dagger)$, and
let $M$ be a finitely generated $S$-module.
Then $M$ is indecomposable as an $R$-module if and only if it is
indecomposable as an $S$-module.
\endproclaim
\demo{Proof}
We know that $M$ is indecomposable as an $R$-module if and only
$\End_R(M)$ has no nontrivial idempotents, and similarly over $S$.
The equality $\End_R(M)=\End_S(M)$ from
Proposition 1.2 now yields the desired result.
\qed
\enddemo

For any ring homomorphism $\varphi: R\to S$, every $S$-module
acquires an $R$-module structure via $\varphi$. We want to
understand when the reverse holds: Given an $R$-module $M$, often
assumed to be finitely generated, when does $M$ have  an $S$-module
structure $(s,x) \mapsto s\circ x$ that is compatible with the
$R$-module structure, that is, $rx=\varphi(r)\circ x$, for $r\in R$
and $x\in M$?    When this happens,
we will say simply that $_RM$ {\it has a compatible } $S$-{\it module structure}.
We are particularly interested in the case where the $S$-module structure is {\it unique}.

\proclaim{Lemma 1.4} Assume $\varphi\colon  R\to S$ satisfies $(\dagger)$.
Let $N$ be a finitely generated $S$-module, and let  $V$ be an
$R$-submodule of $N$.  Then $_RV$ has at most one compatible $S$-module structure.
 In detail:  if $V$ has an $S$-module structure $(s,v)\mapsto s\circ v$
that is compatible with the $R$-module structure on $V$ inherited from the
$S$-module structure $(s,n) \mapsto s\cdot n$ on $N$, then
$s\circ v = s\cdot v$ for all $s\in S$ and $v\in V$.
\endproclaim

\demo{Proof} Let $s\in S$ and $v\in V$ be given.
As before, we fix an integer $t\ge 1$ and choose $r\in R$ such that
$\varphi(r)-s\in \n^t$. Note the following relations:
    $$\n^t\circ V = (\m^t S)\circ V = \m^t\circ(S\circ V)
     = \m^t \circ V = \m^t\cdot V \subseteq \m^t \cdot N =
   \n^t \cdot N
    $$
It follows that we have
$$ s\circ v - s\cdot v = s\circ v - r\circ v + r\cdot v - s\cdot v =
(s-\varphi(r))\circ v + (\varphi(r)-s)\cdot v \in \n^t\circ V + \n^t
\cdot V \subseteq \n^t\cdot  N.
$$
Since $t$ was chosen arbitrarily, we conclude that $s\circ v = s\cdot v$.
\qed\enddemo

\proclaim{1.5. Proposition and Notation}
Assume $\varphi\colon  R\to S$ satisfies $(\dagger)$.
Let $M$ be an $R$-module (not necessarily finitely generated)
that is an $R$-submodule of some finitely generated $S$-module $N$.
Let $\Cal V(M)$ be the set of $R$-submodules of $M$
that have an $S$-module structure compatible with their $R$-module
structure. Then $\Cal V(M)$ is exactly the set of $S$-submodules of
$N$ that are contained in $M$. The set $\Cal V(M)$ has a unique maximal
element ${\V}(M)$.  Moreover, we have ${\V}(M) = \{x\in M \mid Sx\subseteq M\}
= \{x\in N \mid Sx \subseteq M\}$.
\endproclaim
\demo{Proof} The first assertion is clear from (1.4).  It follows that
$\Cal V(M)$ is closed under  sums.  Since $N$ is a Noetherian
$S$-module, the other assertions follow easily. \qed\enddemo

Although $\V(M)$ is defined only when $M$ can be embedded as
an $R$-submodule of some finitely generated $S$-module $N$, its definition is intrinsic.
Thus the submodule $\V(M)$ of $M$ does not depend on the choice of the
module $N$ or the $R$-embedding $M\hookrightarrow N$.  (See Corollary 1.7 for another intrinsic characterization of $\V(M)$.)

\proclaim{Proposition 1.6} Assume $\varphi\colon  R\to S$ satisfies $(\dagger)$, and
let $L$ be an $S$-module (not necessarily
finitely generated). Let $M$ be an $R$-submodule of some finitely
generated $S$-module, and let $\V(M)$ be as in (1.5). Then the
natural injection $\Hom_R(L,\V(M))\to \Hom_R(L,M)$ is an isomorphism.
\endproclaim
\demo{Proof} Let $g\in \Hom_R(L,M)$, and let $W$ be the image of
$f$.  We want to show that $W\subseteq \V(M)$.  Let $h$ be the
composition $L@>g>>M\hookrightarrow N$, where $N$ is some finitely
generated $S$-module containing $M$ as an $R$-submodule. By (1.2),
the map $h$ is $S$-linear, so $W = h(L)$ is an $S$-submodule of $N$.
Therefore we have $W\subseteq \V(M)$.\qed
\enddemo

\proclaim{Corollary 1.7}
Assume $\varphi\colon  R\to S$ satisfies $(\dagger)$.
Let $M$ be an $R$-submodule of a
finitely generated $S$-module.  The following natural maps are
isomorphisms:
$$ \V(M) @>\cong>>\Hom_S(S,\V(M)) @>=>>\Hom_R(S,\V(M))@>\cong>>\Hom_R(S,M)
$$
It follows that $\V(M)$ is exactly the image of the natural map
$\varepsilon\colon\Hom_R(S,M)\to M$ taking $\psi$ to $\psi(1)$.
In particular, if $M$ is finitely generated as an $R$-module, so is $\Hom_R(S,M)$. \qed
\endproclaim

The next result contains the first part of our answer to Question (1) from the Introduction.

\proclaim{Theorem 1.8}  Assume $\varphi\colon  R\to S$ satisfies $(\dagger)$, and let $M$ be
a finitely generated $R$-module.  The following conditions are equivalent:
\roster
\item $M$ has a compatible $S$-module structure.
\item The natural map $\iota\colon  M\to S\otimes_RM$ (taking $x$ to $1\otimes x$) is bijective.
\item The natural map $\varepsilon\colon  \Hom_R(S,M) \to M$ (taking $\psi$ to $\psi(1)$) is bijective.
\endroster
If, in addition, $\varphi$ is  flat, these conditions are equivalent to the following:
\roster
\item"(4)" $S\otimes_RM$ is finitely generated as an $R$-module.
\endroster
\endproclaim
\demo{Proof} The implications (2)$\implies$(1), (3)$\implies$(1), and
(2)$\implies$(4) are clear.
Assume (1), and let $(s,x)\mapsto s\cdot x$ be a compatible $S$-module structure on $M$.  To prove (2), we note that the module $S\otimes_RM$ has two compatible $S$-module structures---the one coming from multiplication in $S$ and the one coming from the $S$-module structure on $M$.  Moreover, with the first structure, $S\otimes_RM$ is finitely generated over $S$.  By (1.4) the two $S$-module structures must be the same.  In particular, for $s\in S$ and $x\in M$ we have $s\otimes x = s(1\otimes x) = 1\otimes (s\cdot x)$.  Therefore the multiplication map $\mu\colon  S\otimes_RM \to M$ (taking $s\otimes x$ to $s\cdot x$) is the inverse of $\iota$.

Still assuming (1), we prove (3).  Since $M$ is finitely generated as an $S$-module,   (1.2) tells us that              $\Hom_R(S,M) = \Hom_S(S,M)$.  Therefore the map $M\to \Hom_S(S,M)$ taking $x\in M$ to the map $s\mapsto s\cdot x$ is the inverse of $\varepsilon$.

(4)$\implies$(2). Assume that $\varphi$ is flat.  By (4), the $S$-module
$S\otimes_RS\otimes_RM$ is finitely generated for the $S$-action on the first variable; therefore its two $S$-module structures (obtained by letting $S$ act on each of the first two factors) are the same, by (1.4).  In particular, $s\otimes t\otimes x = st\otimes 1 \otimes x$ for $s,t\in S$ and $x\in M$.  Therefore the map $S\otimes_R S \otimes_R M \to S\otimes_R M$ taking $s\otimes t\otimes x$ to $st\otimes x$ is the inverse of $1\otimes \iota\colon  S\otimes_R M \to S\otimes_R S \otimes_R M $.  By faithful flatness, $\iota$ is an isomorphism.
\qed\enddemo

In light of Corollary 1.7, we see that the conditions in the previous result are not
equivalent to $\Hom_R(S,M)$ being finitely generated as an $R$-module, even when
$\vf$ is flat.
In the next section, we will show that the ``right" condition is that
 $\Ext^i_R(S,M) $ be finitely generated for  $i=1,\ldots,\dim_R(M)$.

Next we revisit Theorem 1.8 from a slightly different perspective:

\proclaim{Theorem 1.9} Let  $\varphi\colon  R \to S$ be a flat local homomorphism
satisfying $(\dagger)$, and let $M$ be
a finitely generated $S$-module.   The following conditions are equivalent:
\roster
\item $M$ is finitely generated as an $R$-module.
\item The natural map $\iota_M\colon  M\to S\otimes_RM$ (taking $x$ to $1\otimes x$) is bijective.
\item $S\otimes_RM$ is finitely generated as an $R$-module.
\endroster
In particular, if $S$ has a faithful module that is finitely generated as an $R$-module, then $\varphi$ is an isomorphism.
\endproclaim
\demo{Proof}  The implication (1)$\implies$(2) is in Theorem 1.8.  Suppose (2) holds.  The $R$-module $S\otimes_RM$ has two $S$-modules structures and, by (2), is finitely generated with respect to the $S$-action on the second factor.  By Lemma 1.4, the two structures agree, and $S\otimes_RM$ is finitely generated with respect to the $S$-action on the first factor. By faithfully flat descent, $M$ is finitely generated over $R$.  Using (2) again, we get (3).

If (3) holds, then $S\otimes_RM$ is {\it a fortiori} finitely generated for the action of $S$ on the first factor. Again using faithfully flat descent, we get (1).

To prove the last statement, suppose $N$ is a faithful $S$-module that is finitely generated as an
$R$-module.
Let $x_1, \dots, x_t$ generate $N$ as an $S$-module, and define $\alpha:
S\to N^t$ by $1\mapsto (x_1,\dots x_t)$.  The kernel of $\alpha$ is the
intersection of the annihilators of the $x_i$, and this intersection is
$(0)$ since $N$ is faithful.  Thus $S$ embeds in $N^t$ and therefore is
finitely generated as an $R$-module.
Now we put $M = S$ in (2) and note that $\varphi\otimes_RS\colon R\otimes_RS \to S\otimes_R S$ is the composition $R\otimes_R S @>\cong>> S @>\iota_S>>S\otimes_RS$.  Therefore $\varphi\otimes_RS$ is an isomorphism, and by faithful flatness $\varphi$ must be an isomorphism. \qed
\enddemo

\proclaim{Proposition 1.10} Assume $\varphi\colon  R\to S$ satisfies $(\dagger)$.
The following conditions  are equivalent:
\roster
\item $R$ has a compatible $S$-module structure.
\item $\varphi$ is an $R$-split monomorphism.
\item $S$ is a free $R$-module.
\item $\varphi$ is a bijection.
\endroster
\endproclaim
\demo{Proof}
The implication (4)$\implies$(3) is clear.

(1)$\implies$(4).
From (1.8) we conclude that the map $\iota\colon R\to S\otimes_R R$ is bijective,
and it follows that $\varphi$ is the composition of two bijections:
$R @>\iota>> S\otimes_RR \to S$.

(2)$\implies$(1).
Let $\pi\colon S\to R$ be an $R$-homomorphism
such that $\pi\varphi = 1_R$. The composition $\varphi\pi\colon  S\to S$ is $S$-linear by (1.2), so
$\varphi(R)=\varphi\pi(S)$ is an $S$-module, and (1) follows.

(3)$\implies$(2).
Let $B$ be a basis for $S$ as an $R$-module.  Write
$1 = \sum_{i=1}^nr_ib_i$ where the $r_i$ are in $R$ and the $b_i$ are distinct elements of $B$.
If each $r_i$ were in $\m$, we would have $1\in \m S = \n$, contradiction.  Thus we may assume that $r_1$ is a unit of $R$.  Let $\pi:S\to R$ be the $R$-homomorphism taking $b_1$ to $r_1^{-1}$ and $b\in B - \{b_1\}$ to $0$.
Then $\pi\varphi = 1_R$, and we have (2). \qed\enddemo

Now we focus on flat homomorphisms satisfying $(\dagger)$.  (In this context {\it every} finitely generated $R$-module can be embedded in a finitely generated $S$-module, namely $S\otimes_RM$.  Thus $\V(M)$ is always defined.)
Every finite-length $R$-module has a compatible $S$-module structure.  (This follows from (1.12) below, by induction on the length, since $R/\m = S/\m S$.)  There are other examples:

\subhead Example 1.11 \endsubhead Let $R$ be a local ring and $P$ a non-maximal
prime ideal such that $R/P$ is $\m$-adically complete (e.g, $R = (\C[X]_{(X)})[[Y]]$ and  $P = (X)$).  Then $R/P$ has a compatible $\widehat R$-module structure.  Indeed, the map $R/P \to \widehat R/P\widehat R$ is bijective.

\medskip

As we shall see in (1.13), the behavior of prime ideals tells the whole story.  The following lemma is clear from the five-lemma and criterion (2) of (1.8):

\proclaim{Lemma 1.12} Let $\varphi\colon R\to S$ be a flat local
homomorphism satisfying $(\dagger)$, and let
$$
0 \to M' \to M \to M'' \to 0
$$
be an exact sequence of finitely generated $R$-modules.
Then $M$ has a compatible $S$-module structure
if and only if $M'$ and $M''$ have compatible
$S$-module structures.     \qed
\endproclaim

\proclaim{Theorem 1.13} Let $\varphi\colon R\to S$ be a flat local
homomorphism satisfying $(\dagger)$, and
let $M$ be a finitely generated $R$-module.   The following conditions  are equivalent:
 \roster
    \item $M$ has a compatible $S$-module structure.
    \item $S = R+PS$ (equivalently, $R/P$ has a compatible $S$-module structure),
     for every $P\in \Min_R(M)$.
    \item $S = R + PS$ (equivalently, $R/P$ has a compatible $S$-module structure),
    for every  $P\in \Supp_R(M)$.
    \endroster
    \endproclaim
    \demo{Proof} The condition $S = R+PS$ just
    says that the injection $R/P \hookrightarrow
     S\otimes_R(R/P)$ is an isomorphism; now (1.8) justifies the
    parenthetical comments.  If (1)
    holds and $P\in \Min_R(M)$, then there is an injection
    $R/P\hookrightarrow M$, so (1.12)
    with $M' = R/P$ yields (2).  Assume (2).  Given $P\in \Supp_R(M)$
    we have $P\supseteq Q$ for some $Q\in \Min_R(M)$.  Then $R/Q
    \twoheadrightarrow R/P$, and (3) follows from (1.12).
    Assuming (3), choose a prime filtration $M= M_0\subset \dots
    \subset M_t$ with $M_i/M_{i-1} \cong R/P_{i-1}$ with $P_i\in
    \Spec(R)$, $i= 1\dots,t$.  Then $P_i\in \Supp_R(M)$ for each $i$,
    and now (1) follows from (1.12).
    \qed\enddemo

Let  $\varphi\colon (R,\m,k) \hookrightarrow (S,\n,l)$ be a
flat local homomorphism.
Recall that
 $\varphi$ is
{\it separable} if the ``diagonal" morphism $S\otimes_R S\to S$
(taking $a\otimes b$ to $ab$) splits as  $S\otimes_RS$-modules (cf.
\cite{DI}). If, further, $\varphi$ is essentially of finite type,
then $\varphi$ is said to be an {\it \'etale} extension of $R$ (cf.
\cite{Iv}).  An \'etale extension $\vf$ is a  {\it pointed \'etale
neighborhood} of $R$ if $k=l$.  It is easy to see that $\m S = \n$
whenever $\varphi$ is an \'etale extension; thus pointed \'etale
neighborhoods satisfy condition $(\dagger)$.  The $R$-isomorphism
classes of pointed \'etale neighborhoods form a direct system, and
the   Henselization $R\to \Rh$ is the direct limit of them.

\proclaim{Corollary 1.14}
Let $R$ be a local ring
and  $M$  a finitely generated $R$-module.
The following conditions are equivalent.
\roster
\item
$M$ admits an
$\Rh$-module structure that is compatible with its $R$-module structure via
the natural inclusion $R\to \Rh$.
\item
For each $P\in\supp_R(M)$, the ring $R/P$
is Henselian.
\item
For each $P\in\Min_R(M)$, the ring $R/P$
is Henselian.
\item
The ring  $R/\ann_R(M)$ is Henselian. \qed
    \endroster
    \endproclaim

\proclaim{Corollary 1.15}
Let $R$ be a local ring.
The following conditions are equivalent.
\roster
\item
$R$ is Henselian.
\item
For each $P\in\spec(R)$, the ring $R/P$ is Henselian.
\item
For each $P\in\Min(R)$, the ring $R/P$ is Henselian. \qed
    \endroster
    \endproclaim

    \head 2. Vanishing of $\Ext$ \endhead

Our goal in this section is to add a fifth  condition equivalent to the conditions in Theorem 1.8, namely,
that $\Ext^i_R(S,M) = 0$ for $i>0$.  Here $R\to S$ is a flat local homomorphism satisfying $(\dagger)$
and $M$ is a finitely generated $R$-module.
Moreover, we will obtain a sixth equivalent condition, namely, that  $\Ext^i_R(S,M)$ is finitely generated
over $R$ for $i=1,\ldots,\dim_R(M)$.
Since our proof uses complexes we will review the basic yoga here.

\subhead 2.1. Notation and conventions \endsubhead
An {\it $R$-complex} is a sequence of
$R$-module homomorphisms
$$
X = \cdots @>\partial^X_{n+1}>>X_n @>\phantom{i}\partial^X_n\phantom{i}>>X_{n-1}@>\partial^X_{n-1}>>\cdots
$$
such that $ \partial^X_{n-1}\partial^X_{n}=0$ for each integer $n$; the
$n$th {\it homology module} of $X$ is
$\HH_n(X):=\Ker(\partial^X_{n})/\im(\partial^X_{n+1})$.
A complex $X$ is {\it bounded} if $X_n = 0$ for $|n|>>0$, {\it bounded above} if $X_n = 0$ for
$n \gg 0$, and {\it homologically finite} if its total
homology module $\HH(X)=\oplus_n\HH_n(X)$
is a finitely generated $R$-module.

Let $X,Y$ be $R$-complexes.
The  {\it Hom complex} $\Hom_R(X,Y)$ is the $R$-complex defined as
$$\Hom_R(X,Y)_n=\prod_p\Hom_R(X_p,Y_{p+n})$$
with $n$th differential $\partial_n^{\Hom_R(X,Y)}$ given by
$$\{f_p\}\mapsto \{\partial^{Y}_{p+n}f_p-(-1)^nf_{p-1}\partial^X_p\}.$$
A {\it morphism} $X\to Y$
is an element $f = \{f_p\} \in\Hom_R(X,Y)_0$ such that
$\partial^Y_pf_p = f_{p-1}\partial^X_p$ for all $p$, that is,
an element of $\Ker(\partial_0^{\Hom_R(X,Y)})$.

A morphism of complexes $\alpha\colon X\to Y$
induces homomorphisms on homology modules
$\HH_n(\alpha)\colon\HH_n(X)\to\HH_n(Y)$, and $\alpha$ is a
{\it quasi-isomorphism} when each $\HH_n(\alpha)$ is bijective.
The symbol ``$\simeq$'' indicates a quasi-isomorphism.

\subhead 2.2. Base change\endsubhead
Let $\varphi\colon R\to S$ be a flat
homomorphism. For any $R$-complex $X$, the
flatness of $\vf$ provides natural $S$-module isomorphisms
$$\HH_i(S\otimes_R X)\cong S\otimes_R\HH_i(X)$$
for each integer $i$.

\subhead 2.3. A connection with condition $(\dagger)$ \endsubhead
Let $\vf\colon (R,\fm,k)\to (S,\fn,l)$ be a flat local ring homomorphism,
and write $\ol{\vf}\colon k\to S/\fm S$ for the induced ring homomorphism.
Let $X\not\simeq 0$ be an $R$-complex
such that each homology module $\HH_i(X)$ is a finite-dimensional
$k$-vector space, and let $r_i$ denote the vector-space dimension of $\HH_i(X)$.
(In our applications we will consider the case $X=K^R$, the Koszul complex on
a minimal system of generators for $\m$.  By \cite{BH, (1.6.5)}, the
homology $\HH(K^R)$ is annihilated by $\m$, and so
each $\HH_i(K^R)$ is a finite-dimensional $k$-vector space.
Note that $K^R \not \simeq 0$ since $\HH_0(K^R) \cong k$.)
Define $\omega\colon  X\to S\otimes_RX$ by the commutative diagram
$$
\CD
X  @>\omega>>  S\otimes_R X\\
 _{\cong}\searrow  && \nearrow _{\varphi\otimes_RX}\\
& R\otimes_RX
\endCD
$$
where the southeast arrow represents the standard isomorphism.
We have a commutative diagram of
$k$-linear homomorphisms
$$
\CD
\HH_i(X)  @>{\HH_i(\omega)}>> \HH_i(S\otimes_R X)
 @>\cong>> S\otimes_R \HH_i(X) @>\cong>> S\otimes_R k^{(r_i)} \\
@V\cong VV &&&& @V\cong VV\\
k^{(r_i)} & & @>\ol\vf^{(r_i)}>> & & (S/\m S)^{(r_i)}
\endCD
$$
Therefore the morphism $\omega$ is a quasi-isomorphism
if and only if $\ol{\vf}$ is an isomorphism, that is, if and only if
$\varphi\colon  R\to S$ satisfies
the condition $(\dagger)$
of Section 1.

The following result is contained in \cite{FSW07, (5.3)}.

\proclaim{Proposition 2.4} Let $X$ and $Y$ be $R$-complexes such that $\HH_n(X)$ and $\HH_n(Y)$ are finitely generated $R$-modules for each $n$.
Let $\alpha\colon  X\to Y$ be a morphism.  Assume that $P$ is a bounded complex of finitely
generated projective $R$-modules such that $P\not\simeq 0$ and
$\Hom_R(P,\alpha)$ is a quasi-isomorphism.
Then $\alpha$ is a quasi-isomorphism. \qed\endproclaim

We can now put the finishing touch on Theorem 1.8:

\proclaim{Main Theorem 2.5}  Let $\varphi\colon  R\to S$ be a ring
homomorphism satisfying $(\dagger)$, and let $M$ be
a finitely generated $R$-module.  The following conditions are equivalent:
\roster
\item $M$ has a compatible $S$-module structure.
\item The natural map $\iota\colon  M\to S\otimes_RM$ (taking $x$ to $1\otimes x$) is bijective.
\item The natural map $\varepsilon\colon  \Hom_R(S,M) \to M$ (taking $\psi$ to $\psi(1)$) is bijective.
\endroster
If, in addition, $\varphi$ is  flat, these conditions are equivalent to the following:
\roster
\item"(4)" $S\otimes_RM$ is finitely generated as an $R$-module.
\item"(5)" $\Ext^i_R(S,M)$ is a finitely generated $R$-module for $i = 1,\dots,\dim_R(M)$.
\item"(6)" $\Ext^i_R(S,M) = 0$ for all $i > 0$.
\endroster
\endproclaim
\demo{Proof}
The equivalences (1)$\iff$(2)$\iff$(3) are in Theorem 1.8, as is
(3)$\iff$(4) when $\varphi$ is flat.  The implication (6)$\implies$(5) is
trivial, so it remains to assume that $\varphi$ is flat and prove
(5)$\implies$(3)  and (1)$\implies$(6).

(5)$\implies$(3).  Assume that
$\Ext^i_R(S,M)$ is  finitely generated over $R$ for $i = 1,\dots,\dim_R(M)$.
We first show that $\Ext^i_R(S,M)=0$ for each $i>\dim_R(M)$.
Let $P$ be an $R$-projective resolution of $S$,
and set $R'=R/\ann_R(M)$.
The fact that $M$ is an $R'$-module yields the first isomorphism
in the following sequence:
$$
\Hom_R(P,M)\cong\Hom_R(P,\Hom_{R'}(R',M))\cong\Hom_{R'}(P\otimes_R R',M)
\eqno (2.5.1)$$
The second isomorphism is Hom-tensor adjointness.
Of course we have isomorphisms $\HH_n(P\otimes_R R')\cong\Tor_n^R(S,R')$,
so the flatness of $\varphi$ yields
$\HH_n(P\otimes_R R')=0$ for $n > 0$.
Therefore the complex $P\otimes_R R'$ is an $R'$-projective resolution
of $S' :=S\otimes_R R'$.
Since $S'$ is flat over $R'$, we have  $\pd_{R'}(S')\leq\dim(R')$ by a result of Gruson and
Raynaud~\cite{RG, Seconde Partie, Thm.~(3.2.6)}, and
Jensen~\cite{J, Prop.~6}.  Therefore
$\Ext^n_{R'}(S',M)=0$ for each $n> \dim(R') = \dim_R(M)$. This yields the vanishing
in the next sequence, for $n>\dim_R(M)$:
$$
\Ext^n_R(S,M)
\cong \HH_{-n}(\Hom_R(P,M))
\cong \HH_{-n}(\Hom_{R'}(P\otimes_R R',M))
\cong \Ext^n_{R'}(S',M)
=0
$$
The first isomorphism is by definition;
the second one is from (2.5.1);
and the third one is from the fact, already noted, that
$P\otimes_R R'$ is an $R'$-projective resolution
of $S'=S\otimes_R R'$.

Let $I$ be an $R$-injective resolution of $M$.
From Corollary 1.7, it follows that $\Hom_R(S,M)$ is a finitely generated
$R$-module.  Since $\Ext^n_R(S,M)=0$ for $i>\dim_R(M)$ and
$\Ext^n_R(S,M)$ is finitely generated over $R$
for  $1\leq n\leq\dim_R(M)$,  the complex
$\Hom_R(S,I)$ is homologically finite over $R$.

Consider the evaluation morphism $\alpha\colon\Hom_R(S,I)\to I$
given by $f\mapsto f(1)$.  To verify condition (3),
it suffices to show that $\alpha$ is a quasi-isomorphism.  Indeed,
assume for the rest of this paragraph that $\alpha$ is a quasi-isomorphism.
It is straightforward to show that the map
$\HH_0(\alpha)\colon\HH_0(\Hom_R(S,I))\to \HH_0(I)$
is equivalent to the evaluation map
$\varepsilon\colon \Hom_R(S,M)\to M$.  The quasi-isomorphism
assumption implies that $\varepsilon$ is an isomorphism, and so condition (3)
holds.

We now show that $\alpha$ is a quasi-isomorphism.
Let ${\bold x}=x_1,\ldots,x_m$ be a minimal generating sequence for $\m$.
The flatness of $\vf$ conspires with the condition $\m S=\n$ to
imply that $\varphi({\bold x})=\varphi(x_1),\ldots,\varphi(x_m)$ is
a minimal generating sequence for $\n$.
Let $K^R=K^R({\bold x})$ and $K^S=K^S(\varphi({\bold x}))$ denote the
respective Koszul complexes,
and note that we have $\rank_R(K^R_i)=\rank_S(K^S_i)= r:={m\choose i}$.
Let $e_{i,1},\ldots,e_{i,r}$ be an $R$-basis for $K^R_i$,
and let $f_{i,1},\ldots,f_{i,r}$ be
the naturally corresponding $S$-basis for $K^S_i$.
The construction yields a natural isomorphism of $S$-complexes
$\beta\colon K^R\otimes_R S\to K^S$ taking
$e_{i,j}\otimes 1$ to $f_{i,j}$.
On the other hand let $K^{\vf}\colon K^R\to K^S$ be given by
$e_{i,j}\mapsto  f_{i,j}$. By (2.3), the flatness of $\varphi$ and condition
($\dagger$) work together to show that $K^{\vf}$ is a quasi-isomorphism.

The source and target of the morphism $\alpha\colon\Hom_R(S,I)\to I$
 are both homologically finite
$R$-complexes, so  it suffices to verify that the induced morphism
$$\Hom_R(K^R,\alpha)\colon
\Hom_R(K^R,\Hom_R(S,I))\to \Hom_R(K^R,I)$$
is a quasi-isomorphism; see Proposition 2.4.
This isomorphism is verified by the following commutative diagram
$$
\CD
\Hom_R(K^R\otimes_RS,I)  @<\Hom_R(\beta,I)<\simeq< \Hom_R(K^S,I) \\
@V(*)V\cong V  @V\Hom(K^{\vf},I)V\simeq V \\
\Hom_R(K^R,\Hom_R(S,I)) @>\Hom_R(K^R,\alpha)>> \Hom_R(K^R,I)
\endCD
$$
wherein the isomorphism (*) is Hom-tensor adjointness.
The morphism $\Hom_R(\beta,I)$ is a quasi-isomorphism because $I$ is a
bounded-above complex of injective $R$-modules and $\beta$ is a quasi-isomorphism.
(See, e.g., the proof of \cite{Wei, (2.7.6)}.)
The same reasoning shows that $\Hom(K^{\vf},I)$ is a quasi-isomorphism.
From the commutativity of the diagram, it follows that
$\Hom_R(K^R,\alpha)$ is a quasi-isomorphism as well.

(1)$\implies$(6).
Assume that $M$ admits an
$S$-module structure that is compatible with its $R$-module structure via $\vf$.
Since $M$ is finitely generated over $R$, it admits a filtration
by $R$-submodules $0=M_0\subset M_1\subset\cdots\subset M_n=M$
such that $M_i/M_{i-1}\cong R/P_i$ for each $i=1,\ldots,n$
where $P_i\in\supp_R(M)$.  We prove the implication by induction on $n$.

When $n=1$, we have $M\cong R/P$ for some $P\in\spec(R)$.
The implication (1)$\implies$(2)  yields an isomorphism
$M\cong R/P\cong S\otimes_RR/P$.
If $Q$ is an $R$-projective resolution of $S$, the flatness of $S$
implies that $Q\otimes_R R/P$ is an $R/P$-projective resolution
of $S\otimes_RR/P \cong R/P$.
Using the next isomorphisms
$$\Hom_R(Q,R/P)\cong\Hom_R(Q,\Hom_{R/P}(R/P,R/P))
\cong\Hom_{R/P}(Q\otimes_RR/P,R/P)$$
we conclude that
$$\Ext^i_R(S,M)=\Ext^i_R(S,R/P)\cong\Ext^i_{R/P}(R/P,R/P)=0$$
for $i\neq 0$.

Now assume $n>1$ and that the implication holds for each $R$-module $M'$ that admits a
prime filtration with fewer than $n$ inclusions.
Because of the exact sequence
$$0\to M_1\to M\to M/M_1\to 0$$
Lemma 1.12 implies that
$M_1$ and $M/M_1$ admit
$S$-module structures that are compatible with their $R$-module structures via $\vf$.
The induction hypothesis implies
$\Ext^i_R(S,M_1)=0=\Ext^i_R(S,M/M_1)$ for all $i>0$.
Thus, the long exact sequence in $\Ext$ coming from the displayed sequence
implies $\Ext^i_R(S,M)=0$ for all $i>0$.
\qed
\enddemo

\proclaim{Corollary 2.6}
Let $R$ be a local ring and $\fa\subset R$ an ideal.
\roster
\item
The $\fa$-adic completion $\widehat{R}^{\fa}$ is $R$-projective
if and only if $R$ is $\fa$-adically complete.
\item The Henselization $R^h$ is $R$-projective
if and only if $R$ is Henselian.
\item If $R\to R'$ is a pointed {\'e}tale neighborhood
and $R'$ is $R$-projective, then $R = R'$.
\endroster
\endproclaim
\demo{Proof} Suppose  $S:=\widehat{R}^{\fa}$ is $R$-projective.  Putting $M = R$ in Theorem 2.5 and using Proposition 1.10, we see that $R = S$.  This proves (1), and the proofs of (2) and (3) are essentially the same.\qed\enddemo

We conclude this section with several examples showing the necessity of the
hypotheses of Theorem 2.5 with respect to the implications
(5)$\implies$(1) and (6)$\implies$(1).
The examples depend on the following addendum to Proposition 1.10,
in which we no longer assume condition
($\dagger$).

 \proclaim{Proposition 2.7} Let $\varphi\colon A\to B$ be an arbitrary homomorphism of commutative rings.  The following conditions are equivalent:
 \roster
 \item The $A$-module $A$ has a $B$-module structure $(b,a)\mapsto b\circ a$ such that
 $$a_1a_2 = \varphi(a_1)\circ a_2 \ \text{for all}\  a_1, a_2 \in A.\tag i$$
 \item $A$ is a ring retract of $B$, that is, there is a ring homomorphism $\psi\colon B\to A$ such that
 $$\psi\varphi(a) = a \ \text{for each}\ a\in A.\tag ii$$
 \endroster
 These conditions imply that $\varphi$ is an $A$-split injection.
 \endproclaim
 \demo{Proof} Assuming (1), we define a function $\psi\colon  B \to A$ by by putting $\psi(b) := b\circ1_A$ for each $b\in B$.  Condition (ii) follows immediately from (i).  Also, given $a\in A$ and $b\in B$ we have
 $\psi(ab)= \psi(\varphi(a)b) = (\varphi(a)b)\circ1_A = \varphi(a)\circ (b\circ1_A)$, by associativity of the $B$-module structure.  Condition (i)
 implies $\varphi(a)\circ (b\circ1_A)   = a(b\circ1_A) = a\psi(b)$,  so $\psi$ is $A$-linear.   This shows that $\varphi$ is an $A$-split injection.

 Still assuming (1), let $b_1,b_2\in B$.   By associativity of the $B$-module structure, we have
 $$\psi(b_1b_2) =  (b_1b_2)\circ1_A = b_1\circ(b_2\circ1_A) = b_1\circ\psi(b_2).\tag iii$$ On the other hand, the $A$-linearity of $\psi$ yields
$\psi(b_1)\psi(b_2) = \psi(b_1\varphi(\psi(b_2)))$.  By  (iii),
this implies $\psi(b_1\varphi(\psi(b_2))) = b_1\circ\psi\varphi\psi(b_2) = b_1\circ \psi(b_2)$.
Thus $\psi(b_1)\psi(b_2) = b_1\circ\psi(b_2)$, and so (iii)
implies that $\psi$ is a ring homomorphism.

 For the converse, assume (2), and set $b\circ a := \psi(b\varphi(a))$ for all $a\in A$ and $b\in B$.  One checks readily the equalities
 $(b_1b_2)\circ a = a\psi(b_1b_2) = b_1\circ(b_2\circ a)$ for $b_i\in B$ and $a\in A$.
Thus we have defined a legitimate $B$-module structure on $A$.
The verification of (i) is easy and left to the reader.\qed
 \enddemo

Our first example shows why
we need to assume that the induced map between the residue
fields of $R$ and $S$ is an isomorphism
in the implications
(5)$\implies$(1) and (6)$\implies$(1) of Theorem 2.5.

\subhead Example 2.8 \endsubhead
Let
$\vf\colon K\to L$ be a proper field extension.
Then $\vf$ is a flat local homomorphism
and  $\fm_K L=\fm_L$ (but the induced map $K/\fm_K\to L/\fm_L$ is
not an isomorphism).  If we take $M$ = $R$, then conditions (5) and (6) of Theorem 2.5 are satisfied, but (1) is not.  Indeed, suppose (1) holds.   Proposition 2.7 provides a field homomorphism $\psi\colon L\to K$ such that $\psi\varphi $ is the identity map on $K$.  Since $\psi$ is necessarily injective, it follows that $\psi$ and $\varphi$ are reciprocal isomorphisms, contradiction.\qed
\medskip
The next example shows the necessity of the condition
$\fm S= \fn$ for the implications (5) $\implies$ (1) and (6)$\implies$(1) in Theorem 2.5.

\subhead Example 2.9 \endsubhead
Let $k$ be a field and $p \ge 2$ an integer. Set $R = k[[X^p]]$ and $S=k[[X]]$, and let $\varphi\colon R\to S$
be the inclusion map. Again, we put $M= R$.  Then $\varphi$ is a local homomorphism inducing an isomorphism on residue fields (but $\m_RS\ne \m_S$).
Since $S$ is a free $R$-module (with basis $\{1,X,\dots,X^{p-1}\}$),   conditions (5) and (6)
are satisfied.  Suppose, by way of contradiction, that (1) is satisfied. Using Proposition 2.7, we get a ring homomorphism $\psi\colon S\to R$ such that $\psi\varphi$ is the identity map on $R$.    Putting $z := \psi(X)$, we see that $X^p = \psi(z^p) \in \m_R^p$, an obvious contradiction.

Similarly, let $R$ be a regular local ring of characteristic $p>0$.
Take $S=R$ and assume that $R$ is F-finite, that is, that
the Frobenius endomorphism $\varphi\colon R\to S$
is module-finite.  (This holds, for example, if $R$ is a power series
ring over a perfect field.)
As an $R$-module, $S$ is  flat by \cite{K}, and therefore
free.  Thus conditions (5) and (6) hold.  Assume $k:=R/\m_R$ is perfect and that
$\dim(R) > 0$.  Then $\varphi$ induces an isomorphism on residue fields, and
essentially the same argument as above shows that condition (1) fails.

\medskip

The next two examples show why we need $\vf$ to be flat for the implications
(5) $\implies$ (1) and (6) $\implies$ (1), respectively.  Note that the homomorphism $\vf$ satisfies ($\dagger$) in both examples and has
 finite flat dimension in Example 2.10.

\subhead Example 2.10 \endsubhead
Let $R$ be a local ring with $\depth(R)\geq 1$ and fix an $R$-regular element
$x\in\fm$.  We consider the natural surjection $\vf\colon R\to R/(x)$.
It is straightforward to show $\Ext^1_R(R/(x),R)\cong R/(x)$ and
$\Ext^n_R(R/(x),R)=0$ when $n\neq 1$.  In particular, each $\Ext^n_R(R/(x),R)$
is finitely generated over $R$.  Suppose (1) holds, and let $\psi\colon S\to R$ be the retraction promised by Proposition 2.7.  Then $x = \psi\varphi(x) = 0$, contradition.

\medskip

\subhead Example 2.11 \endsubhead
Let $R$ be a local Artinian Gorenstein ring with residue field $k\neq R$.
We consider the natural surjection $\vf\colon R\to k$.
Because $R$ is self-injective, we have $\Ext^n_R(k,R)=0$ when $n\neq 0$.
Thus conditions (5) and (6) of Theorem 2.5 hold.  As in Example 2.10, we see easily that (1) fails.

\medskip

\head 3. Extended modules \endhead

Let $\varphi\colon (R,\m) \to (S,\n)$ be a flat local homomorphism.
Given a finitely generated $S$-module $N$, we say that $N$ is {\it extended} provided there is an $R$-module $M$ such that $S\otimes_RM \cong N$ as $S$-modules.  By faithfully flat descent, such a module $M$, if it exists, is unique up to $R$-isomorphism and is necessarily finitely generated.

We begin with a ``two-out-of-three" principle, which is well known when
$S=\widehat R$.
The proof in general seems to require a different approach from the proof in that special case.  The following notation will be used in the proof:  Given a ring $A$ and $A$-modules $U$ and $V$, we write $U \mid_A V$ to indicate that $U$ is isomorphic to a direct summand of $V$.

\proclaim{Proposition 3.1} Let $\varphi\colon R \to S$ be a flat local homomorphism.  Let $N_1$ and $N_2$ be finitely generated $S$-modules, and put $N = N_1\oplus N_2$.  If two of the modules
$N_1, N_2, N$ are extended, so is the third.
\endproclaim
\demo{Proof}  We begin with a claim: If $M_1$ and $M$ are finitely generated $R$-modules, and if $S\otimes_RM_1\mid_S S\otimes_RM$, then $M_1\mid_RM$.  To prove the claim, write $S\otimes_RM \cong (S\otimes_RM_1)\oplus U$. We assume, temporarily, that $R$ is Artinian.  By \cite{Wi98, (1.2)} we know, at least, that there is some positive integer $r$ such that $M_1\mid_R M^{(r)}$ (a suitable direct sum of copies of $M$).
Write $M_1\cong \oplus_{i=1}^s V_i$ where each $V_i$ is indecomposable.  We proceed by induction on $s$.  Since $V_1\mid_R M^{(r)}$,  the Krull-Remak-Schmidt theorem (for finite-length modules) implies that $V_1\mid_RM$, say, $V_1 \oplus W \cong M$.
This takes care of the base case $s=1$.  For the inductive step, assume $s>1$ and set
$W_1=\oplus_{i=2}^s V_i$.
We have $V_1\oplus W \cong M$ and $M_1\cong V_1\oplus W_1$, and hence
$$(S\otimes_R V_1) \oplus (S\otimes_R W) \cong S\otimes_RM \cong (S\otimes_RM_1)\oplus U \cong (S\otimes_RV_1)\oplus (S\otimes_RW_1) \oplus U.$$
Direct-sum cancellation \cite{Ev} implies
$(S\otimes_R W) \cong (S\otimes_RW_1) \oplus U$.
The inductive hypothesis, applied to the pair $W_1,W$, now implies that $W_1\mid_R W$; therefore $M_1\mid_RM$.  This completes the proof of the claim when $R$ is Artinian.

In the general case, let $t$ be an arbitrary positive integer, and consider the flat local homomorphism $R/\m^t \to S/\m^t S$.  By the Artinian case, $M_1/\m^t M_1 \mid_{R/\m^t}M/\m^t M$.  Now we apply Corollary 2 of \cite{G} to conclude that $M_1 \mid_R M$, as desired.

Having proved our claim, we now complete the proof of the proposition.  If $N_1$ and $N_2$ are extended, clearly $N$ is extended.  Assuming $N_1$ and $N$ are extended, we will prove that $N_2$ is extended.  (The third possibility will then follow by symmetry.)  Let $N_1 \cong S\otimes_RM_1$ and $N\cong S\otimes_RM$.    Thus $S\otimes_RM_1 \mid_S S\otimes_RM$, and by the claim there is an $R$-module $M_2$ such that $M_1\oplus M_2 \cong M$. Now
$N_1\oplus (S\otimes_RM_2) \cong S\otimes_RM \cong N_1\oplus N_2$, and, by direct-sum cancellation \cite{Ev}, we have $S\otimes_RM_2 \cong N_2$.
\qed\enddemo



There is a ``two-out-of-three" principle for short exact sequences as well, though some restrictions apply.  Variations on this theme have been used in the literature, e.g., in \cite{CPST}, \cite{LO}, \cite{Wes}.

\proclaim{Proposition 3.2} Let $\vf\colon R \to S$ be a flat local
homomorphism satisfying $(\dagger)$, and consider  an exact
sequence of finitely generated $S$-modules $0 \to N_1 \to N \to N_2 \to 0$. \roster
    \item Assume that $N_1$ and $N_2$ are extended. If $\Ext^1_S(N_2,N_1)$
    is finitely generated as an $R$-module,
    then $N$ is extended.
    \item Assume that $N$ and $N_2$ are extended. If $\Hom_S(N,N_2)$ is finitely generated
    as an $R$-module, then $N_1$ is extended.
        \item Assume that $N_1$ and $N$ are extended. If $\Hom_S(N_1,N)$ finitely generated as an $R$-module, then $N_2$ is
    extended.
    \endroster
\endproclaim
    \demo{Proof} For (1), let $N_i = S\otimes_RM_i$  where the $M_i$ are finitely generated
    $R$-modules.  We have natural homomorphisms $\Ext^1_R(M_2,M_1) @>\alpha>>S\otimes_R\Ext^1_R(M_2,M_1)
@>\beta>> \Ext^1_S(N_2,N_1)$.  The map $\beta$ is an isomorphism because $\vf$ is flat,
    $M_2$ is finitely generated and $R$ is Noetherian.  Therefore $S\otimes_R\Ext^1_R(M_2,M_1)$ is finitely generated as an $R$-module, and now Theorem 1.8 ((4)$\implies$(2)) says that $\alpha$ is an isomorphism.  This means that the given exact sequence
    of $S$-modules is isomorphic to $S\otimes_R{\bold M}$
    for some exact sequence of $R$-modules ${\bold M}=(0 \to M_1 \to M \to M_2
    \to 0)$.  Clearly, this implies $S\otimes_RM \cong N$.

    To prove (2), let $N\cong S\otimes_RM$ and $N_2 \cong S\otimes_RM_2$, where $M$ and $M_2$ are finitely generated $R$-modules.  Essentially the same proof as in (1), but with $\Hom$ in place of $\Ext$, shows that the given homomorphism $N\to N_2$ comes from a homomorphism $f\colon  M\to M_2$.  Then $M_1 \cong S\otimes_R\Ker(f)$.

For (3), we let $N_1 \cong S\otimes_RM_1$ and $N\cong S\otimes_RM$; we deduce that the given homomorphism $N_1 \to N$ comes from a homomorphism
 $g\colon M_1\to M$.  Then $N_2 \cong S\otimes_R\Coker(g)$.\qed
    \enddemo

Here is a simple application  of Part (1) of Proposition 3.2 (cf. \cite{LO} and \cite{Wi01} for much more general results):

  \proclaim{Proposition 3.3} Let $(R,\m)$ be a one-dimensional local ring whose $\m$-adic completion $S = \widehat R$ is a domain.  Then every finitely generated $S$-module is extended.
  \endproclaim
  \demo{Proof} Given a finitely generated $S$-module $N$, let $\{x_1,\dots, x_n\}$ be a maximal
  $S$-linearly independent subset of $N$.  The submodule $F$ generated by the $x_i$ is free and therefore extended.  The quotient module $N/F$ is torsion and hence of finite length.  Therefore $N/F$ is extended.  Since $\Ext^1_S(N/F,F)$ has finite length,
the module $N$ is extended, by (3.2).\qed\enddemo

Notice that Part (1) of Proposition 3.2 applies also when $N_2$ is free on the punctured spectrum.
For in this case $\Ext^1_S(N_2,N_1)$ has finite length over $S$ and therefore is finitely generated
as an $R$-module.  A more subtle condition  that forces  $\Ext^1_S(N_2,N_1)$ to have finite length is
that there are only finitely many isomorphism classes of modules $X$ fitting into a short exact
sequence $0\to N_1\to X\to N_2 \to 0$.  (Cf. \cite{CPST, (4.1)}.)

\medskip

Of course not every module over the completion, or over the Henselization, is extended.  Suppose, for example, that $R= \C[X,Y]_{(X,Y)}/(Y^2-X^3-X^2)$, the local ring of a node.
Then $R$ is a domain, but $\widehat R \cong \C[[U,V]]/(UV))$, which has two minimal prime ideals
$P = (U)$ and $Q = (V)$.  Since $R$ is a domain, any extended $\widehat R$-module $N$ must have the property that $N_P$ and $N_Q$ have the same vector-space dimension (over $\widehat R_P$ and $\widehat R_Q$, respectively).  Thus the $\widehat R$-module $\widehat R/P$ is not extended.  (This behavior was the basis for the first example of failure of the Krull-Remak-Schmidt theorem for finitely generated modules over local rings. See the example due to R.\ G.\ Swan in \cite{Ev}.  The idea is developed further in \cite{Wi01}.) The module $\widehat R/P$ is free on the punctured spectrum and therefore, by Elkik's theorem \cite{El}, is extended from the Henselization $\Rh$.  With $\widehat R/P \cong \widehat R\otimes_{\Rh}V$, we see that the $\Rh$-module $V$ is not extended from $R$.

\medskip

Next, we turn to the question of whether every finitely generated module over $S$ is a direct summand of a finitely generated extended module.  This weaker property is often useful in questions concerning ascent of finite representation type (cf. \cite{Wi98, Lemma 2.1}). Although the next result is not explicitly stated in \cite{Wi98}, the main ideas of the proof occur there.
Note that we do not require that $R/\m = S/\n$.

\proclaim{Theorem 3.4} Let $\vf\colon R\to S$ be a flat local homomorphism,
and assume $S$ is separable over $R$ (that is, the diagonal map $S\otimes_RS \to S$ splits as $S\otimes_RS$-modules).  Then every finitely generated $S$-module is a
direct summand of a finitely generated extended module.
\endproclaim
\demo{Proof} Given a finitely
generated $S$-module $N$, we apply $-\otimes_S N$
to the diagonal map, getting a split surjection of $S$-modules
$\pi\colon S\otimes_RN \twoheadrightarrow N$, where the $S$-module structure on $S\otimes_RN$
comes from the $S$-action on $S$, not on $N$.  Thus we have a split injection of $S$-modules
 $j\colon N\to S\otimes_RN$.  Now write $N$ as a direct
union of finitely generated $R$-modules $M_i$.  The flatness of
$\vf$ implies that $S\otimes_R N$ is a direct union of the modules $S\otimes_R
M_i$.
The finitely generated $S$-module $j(N)$ must be
contained in some $S\otimes_R M_i$. Since $j(N)$ is a direct summand
of $S\otimes_R N$, it must be a direct summand of the smaller
module $S\otimes_RM_i$.\qed
\enddemo

\proclaim{Corollary 3.5}
Let $R\to \Rh$ be the Henselization of the local ring $R$.  Then every finitely generated $\Rh$-module is a direct summand of a finitely generated extended module.
\endproclaim
\demo{Proof} Let $N$ be a finitely generated $\Rh$-module.  Since $R\to \Rh$ is a direct limit of  \'etale neighborhoods $R \to S_i$, $N$ is extended from some $S_i$. Now apply Theorem 3.4 to the extension $R\to S_i$. \qed
\enddemo

The analogous result can fail for the completion:

\subhead Example 3.6 \endsubhead  Let $(R,\m)$ be a countable local ring of
dimension at least two.  Then $R$ has only countably many
isomorphism classes of finitely generated modules. Using the Krull-Remak-Schmidt
theorem over the $\m$-adic completion $\widehat R$, we see that
only countably many isomorphism classes of indecomposable $\widehat R$-modules occur in
direct-sum decompositions of finitely generated
extended modules. We claim, on the other hand, that $\widehat R$ has
uncountably many isomorphism classes of finitely generated
indecomposable modules.  To see this, we recall that $\widehat R$,
being complete, has countable prime avoidance; see \cite{SV}.
By Krull's principal
ideal theorem, the maximal ideal of $\widehat R$ is the union of the
height-one prime ideals.  It follows that $\widehat R$ must have
uncountably many height-one primes $P$, and the $\widehat R$-modules
$\widehat R/P$ are pairwise non-isomorphic and indecomposable.

\medskip

If $\varphi\colon (R,\m,k) \to (S,\n,l)$ is flat and satisfies $(\dagger$), we know that every finite-length $S$-module is extended.  We close with an example showing that the condition $k=l$  cannot be deleted, even for a module-finite  \'etale  extension of Artinian local rings.

\subhead Example 3.7\endsubhead Let $R = \R[X,Y]/(X,Y)^2$ and $S = \C\otimes_{\R}R = \C[X,Y]/(X,Y)^2$. We claim that, for $c\in \C$, the module $N:=S/(X+cY)$ is extended (if and) only if $c\in \R$.   The minimal presentation of $N$ is $S@>X+cY>>S\to N \to 0$.  If $N$ were extended, the $1\times 1$ matrix $X+cY$ would be equivalent to a matrix over $R$.  In other words, we would have $X+cY = u(r+sX+tY)$ for some unit $u$ of $S$ and suitable elements $r,s,t\in \R$.  Writing $u=a+bX+dY$, with $a,b,d\in \C$ and $a\ne 0$, we see, by comparing coefficients of $1,X$ and $Y$, that $c = t/s \in \R$.

\enddocument